\documentclass[12pt]{article}
\begin{document}
\def\C{\mathbf{C}}
\def\R{\mathbf{R}}
\def\bC{\mathbf{\overline{C}}}
\def\Z{\mathbf{Z}}
\def\T{\mathbf{T}}
\def\Ima{\mathrm{Im}\, }
\title{Metrics of constant positive curvature with conic singularities. A survey}
\author{Alexandre Eremenko}
\maketitle
\begin{abstract} We consider conformal metrics of constant curvature $1$
on a Riemann surface, with finitely many prescribed conic singularities
and prescribed angles at these singularities. 
Especially interesting case which was studied by C. L. Chai, C. S Lin and
C. L. Wang is described in some detail, with simplified proofs.


\end{abstract}

\noindent
{\bf 1. Introduction.}
\vspace{.2in}

This is an expanded version of the talk given by the author
on the conference Geometry, Differential Equations and Analysis,
in Kharkiv, June 17, 2019.

Let $S$ be a compact Riemann surface and $A=\{ a_1,\ldots,a_n\}\subset S$
a finite set. A conformal Riemannian metric is given by the length element
$\rho(z)|dz|,\;\rho>0$, and the curvature of this metric is
$$\kappa=-\frac{\Delta\log\rho}{\rho^2}.$$
We assume that the curvature is constant on $S\backslash A$ while at the
points $a_j$ the metric has {\em conic singularities}:
$$\rho(z)\sim c|z-a_j|^{\alpha_j-1},\quad z\to a_j.$$
Positive numbers $\alpha_j$ are the {\em angles} at the singularities
({\em we measure angles in turns; one turn $=2\pi$ radians}). 

When there are no singularities, a classical result says that
on each Riemann surface there is a metric with $\kappa=1$
when $S$ is the sphere, $\kappa=0$ when $S$ is a torus, and $\kappa=-1$
for all other Riemann surfaces. This metric is unique except for $\kappa=0$
when it can be multiplied by an arbitrary positive number.

The problem we discuss here is how to understand (describe, classify)
such metrics with prescribed singularities and angles at the singularities.
This problem has been studied by a large variety of methods:
PDE \cite{P1,P2,P3,P4,LW,LT,T1}, non-linear functional analysis \cite{BdMM,CL},
analytic theory of linear ODE \cite{E1,ER,ET,L}, complex analysis
\cite{E2,EGT1} geometry
\cite{EGT1,MP,MP1,LT}, elliptic functions and modular forms \cite{EG,LW},
and holomorphic dynamics \cite{BE}.

The simplest examples of such metrics are obtained from polygons.
An $n$-gon is a simply connected bordered Riemann surface equipped
with a Riemannian metric of constant curvature $\kappa$, and such that
the boundary consists of $n$ geodesic arcs meeting at $n$ corners.
Two $n$-gons are called conformally equivalent if there is a conformal
homeomorphism between them sending vertices to vertices. Gluing
a  polygon
to its mirror image we obtain a sphere with metric of constant curvature
with conic singularities at the corners. When $\kappa=0$, the Gauss--Bonnet
theorem implies that the sum of the interior angles at the corners
is $n-2$, and  the
Christoffel--Schwarz formula tells us thats for prescribed angles
and prescribed conformal class, there is one polygon,
up to similarity. The situation is similar in hyperbolic geometry,
but it is much more complicated in spherical geometry: there
are additional restrictions on the angles, and when they are satisfied
the polygon with prescribed angles in a given conformal class
may be not unique. Metrics on the sphere coming from polygons
are characterized by the property that there is a circle
on the sphere which contains all singularities, and the metric
is symmetric with respect to this circle. 

A complete solution of our problem is known when $\kappa\leq 0$.
The Gauss--Bonnet theorem implies that on a surface of Euler characteristic
$\chi(S)$,
\begin{equation}\label{gb}
\chi(S)+\sum_{j=1}^n(\alpha_j-1)=\frac{1}{2\pi\,}(\mbox{total integral curvature}),
\end{equation}
{\em so the expression in the LHS
has the same sign as $\kappa$}. 

When $\kappa\leq 0$, this is a necessary and sufficient condition for the
existence of the metric; it is unique when $\kappa<0$ and unique up
to a constant factor when $\kappa=0$.

This result goes back to Picard, \cite{P1,P2,P3,P4}. Modern proofs can be
found in \cite{H,T1}. All these proofs 
are based on the consideration
of solutions of the differential equation
$$\Delta u+\kappa e^{2u}=0$$
on $S\backslash A$ with prescribed behavior at the singularities.
Here $u=\log\rho$.

On the other hand, the problem with $\kappa>0$ is wide open,
and it is the subject of this paper.
\vspace{.1in}

{\em From now on we assume that $\kappa=1.$}
\vspace{.2in}

\noindent
{\bf 2. Developing map and a general conjecture.}
\vspace{.2in}

A surface of curvature $1$ is locally isometric to a  region
on the standard sphere which we denote by $\bC$. (The length element of
the metric on $\bC$ is $2|dz|/(1+|z|^2)$.) This local isometry
can be analytically continued along any path not passing through
the singularities so we have a multivalued {\em developing map} 
$$f:S\backslash A\to \bC.$$
As a local isometry, this map is holomorphic,
and its monodromy is a subgroup of the group of isometries
(rotations) of the sphere $SO(3)=PSU(2)=SU(2)/\{\pm I\}$.
Near the singularities we have
\begin{equation}\label{s}
f(z)\sim c(z-a_j)^{\alpha_j}.
\end{equation}
Conversely, any multivalued locally biholomorphic
function $f$ on $S\backslash A$ with
$PSU(2)$ monodromy and satisfying (\ref{s}) near the points $a_j\in A$
is a developing map of some metric on the sphere of curvature $1$ with
conic singularities at $a_j$ with angles $\alpha_j$. The metric
is recovered from $f$ by the formula
$$\rho(z)=\frac{2|f'|}{(1+|f|^2)}.$$
Developing map is not unique: two developing maps $f$ and $g$
define the same metric if $f=\phi(g)$ where $\phi\in PSU(2)$.

So classification of metrics is equivalent to classification of
developing maps modulo composition with rotations.
 
For a developing map $f$ it can happen that $\phi(f)$ with $\phi\in PSL(2,\C)$
is also a developing map even when $\phi$ is not a rotation.
This happens if and only if $\phi$ conjugates the monodromy group $\Gamma$
of $f$
to a subgroup of $PSU(2)$. It is easy to see that this can happen
with $\phi\not\in SU(2)$ if and only
if 
$\Gamma$ is isomorphic to a subgroup of the unit circle $O(2)$.

In this case, the monodromy and the metric are
called {\em co-axial}.
Existence of co-axial metrics justifies the following definition:
\vspace{.1in}

{\em Two metrics are called equivalent if their developing maps
are obtained from each other by post-composition with a linear-fractional
transformation.} 
\vspace{.1in}

Equivalence classes are $3$-parametric when $\Gamma$ is trivial, 
one parametric when $\Gamma$ is a non-trivial subgroup
of the unit circle, and 
consist of one element in all other cases.

Now we can state a general conjecture.
\vspace{.1in}

{\bf Conjecture 1.} {\em For any compact Riemann surface $S$,
and any prescribed singularities $a_j$ and angles $\alpha_j$,
there are finitely many equivalence classes
of metrics of curvature $1$
on $S$ with these angles at these singularities.}
\vspace{.1in}

Even in the simplest cases there can be more than one class, in sharp contrast
with the case of non-positive curvature. The simplest example
of non-uniqueness
is given in Section 6. Conjecture 1 has been proved for the case when
$S$ is the sphere with $4$ singularities \cite{ER}.

If the general conjecture is true, the next question is 
\vspace{.1in}

{\bf Question 1.} {\em How many equivalence classes of metrics exist on a compact Riemann
surface
with prescribed angles and singularities?}
\vspace{.1in}

The proof in \cite{ER} is not constructive and gives no upper estimate.
\vspace{.2in}

\noindent
{\bf 3. General restrictions on the angles.}
\vspace{.2in}

Here we address the question what angles
$\alpha_j$ can occur (without prescribing the position of the
singularities $a_j$ or conformal type of $S$).
The Gauss-Bonnet theorem for the sphere gives
\begin{equation}\label{gb1}
2+\sum_{j=1}^n(\alpha_j-1)>0.
\end{equation}
Unlike in the case $\kappa\leq 0$, in the case when $S$ is the sphere,
there is another condition, which is
called the {\em closure condition}:
\begin{equation}\label{cc}
d_1({\bf{\alpha-1}},\Z_o^n)\geq 1,
\end{equation}
where ${\bf{\alpha-1}}=(\alpha_1-1,\ldots,\alpha_n-1)$, $\Z_o^n$ is the odd
integer lattice (the set of points in $\R^n$ whose coordinates are integers
whose sum is odd), and $d_1$ is the $\ell_1$ distance.
\vspace{.1in}

\noindent
{\bf Theorem 1.} (Mondello and Panov \cite{MP}) {\em
Conditions (\ref{gb1}) and (\ref{cc}) are necessary for existence of
a metric of curvature $1$ on the sphere with angles $\alpha_j$
and some (unspecified)
singularities
$a_j$.

Conditions (\ref{gb1}) and (\ref{cc}) with strict inequality are
sufficient.

Equality in condition (\ref{cc}) can only hold
for metrics with co-axial monodromy.}
\vspace{.1in}

Several special cases of (\ref{cc}) were known before \cite{MP},
usually they were stated in different forms.

Possible angles of co-axial metrics on the sphere are described in \cite{E2}.

To state the result, let us call a vector
${\bf{\alpha}}=(\alpha_1,\ldots,\alpha_n)$ with positive coordinates admissible
if there exists a co-axial metric with angles $\alpha_j$,
and suppose without loss of generality that $\alpha_1,\ldots,\alpha_m$ are not integers,
while $\alpha_{m+1},\ldots,\alpha_n$ are integers.
\vspace{.1in}

\noindent
{\bf Theorem 2.} {\em For ${\bf{\alpha}}$ to be admissible it is necessary
that:
\vspace{.1in}

\noindent
(i) there exist a choice of signs $\epsilon_j\in\{\pm1\}$ and
an integer $k'\geq 0$ such that
$$\sum_{j=1}^m\epsilon_j\alpha_j=k',$$

\noindent
(ii) the integer
$$k'':=\sum_{j=m+1}^n\alpha_j-n-k'+2\quad\mbox{is even and non-negative}.$$

If the numbers 
$${\bf{c}}=(c_1,\ldots,c_q):=(\alpha_1,\ldots,\alpha_m,
\underbrace{1,\ldots,1}_{ k'+k''\;\mathrm{times}})$$
are incommensurable, then (i) and (ii) are also sufficient.
\vspace{.1in}

\noindent
(iii) If ${\bf{c}}=\eta{\bf{b}}$, where $\eta\neq 0$ and $b_j$ are
integers, then there is an additional necessary condition
\begin{equation}\label{ac}
2\max_{m+1\leq j\leq n}\alpha_j\leq\sum_{j=1}^q|b_j|,
\end{equation}
and in this case the three conditions (i), (ii) and (\ref{ac})
are sufficient.}
\vspace{.1in}

Parameter count shows that in the case of co-axial monodromy, if the
number of non-integer angles $m$ is greater than $2$, the positions
of the singularities cannot be prescribed.

These results give a complete description of possible angles at
the conic singularities for metrics of curvature $1$ on the sphere.

The situation on surfaces of higher genus is simpler:
\vspace{.1in}

\noindent
{\bf Theorem 3.} (Mondello and Panov \cite{MP1}) {\em For any even $\chi\leq 0$
and any $\alpha_j,\; j=1,\ldots,n$ satisfying
$$\chi+\sum_{j=1}^n(\alpha_j-1)>0$$
there exists a compact Riemann surface $S$ of Euler characteristic $\chi$
and a metric of curvature $1$ on $S$ with conic singularities
with angles $\alpha_j$. Moreover, this metric has non-coaxial
monodromy unless $\chi=0$ and all $\alpha_j$ are integers.}
\vspace{.2in}

\noindent
{\bf 4. Fuchsian differential equations.}
\vspace{.2in}

The monodromy of the developing map
is a subgroup of the group of linear fractional transformations $PSL(2,\C)$,
therefore the Schwarzian derivative
$$F:=\frac{f'''}{f}-\frac{3}{2}\left(\frac{f''}{f}\right)^2$$
is single-valued and 
defines a quadratic differential $F(z)dz^2$ holomorphic on $S\backslash A$,
having double poles at the singularities.

We conclude that $f=w_1/w_2$, where
$w_1$ and $w_2$ are linearly independent solutions of a differential equation
\begin{equation}\label{de}
w''+Pw'+Qw=0,
\end{equation}
where
$$F=-P'-P^2/2+2Q,$$
and there is
some freedom in this choice. Changing $P$ results in
multiplication of $w_1$ and $w_2$ by a common factor.
For example we can choose $P=0, Q=F/2$,
so that (\ref{de})
becomes
\begin{equation}\label{de1}
w''+(F/2)w=0.
\end{equation}

From the asymptotics of $f$ at the points $a_j$ we conclude that
equation (\ref{de}) has regular singularities, 
and the exponent differences at each singularity $a_j$ is $\pm\alpha_j$.
Moreover, the projective monodromy of this equation must be conjugate
to a subgroup
of $PSU(2)$.

We have a bijective correspondence between equivalence classes
of metrics of curvature $1$
with conic singularities at $a_j$ with angles $\alpha_j$ and
differential equations (\ref{de1}) with singularities at $a_j$,
exponent differences $\pm\alpha_j$ and $PSU(2)$ projective monodromy
up to conjugacy.

The set of equations (\ref{de1}) on a given Riemann surface of genus $g$
with prescribed singularities and exponent
differences depends on $3g+n-3$ parameters which are called {\em accessory
parameters}. These parameters have to be chosen so that the projective
monodromy
is conjugate to a subgroup of $PSU(2)$.
So our main question is equivalent to the following
\vspace{.1in}

{\em For equation (\ref{de1}) with prescribed singularities
and prescribed real exponent differences, how many choices of
accessory parameters exist so that the monodromy of this equation
is conjugate to a subgroup of $PSU(2)$?}
\vspace{.1in}

In the simplest cases of the sphere with four singularities and a torus
with one singularity
we have one (complex) equation on one accessory
parameter. 

Below, in sections 6--10,
we describe all cases when the answer to the main questions stated in
Section 3 is known.
But first we state some general results.
\vspace{.1in}

\noindent
{\bf Theorem 4.} (Feng Luo, \cite{L}).
{\em Let $Q$ be the fibration over the Teichm\"uller space $T_{g,n}$
of surfaces
$S$ of genus $g$ with $n$ punctures whose fiber at a point $S$
is the space of
quadratic differentials with at most double poles at the punctures.
Let
$$p:Q\to{\mathrm{Hom}}\left(\pi_1(S),PSL(2,\C)\right)/PSL(2,\C)$$
be the monodromy map.  Then $p$ is locally biholomorphic at every
point $x$ where all the exponent differences are not integers, and
$p(x)$ is a smooth point
of ${\mathrm{Hom}}\left(\pi_1(S),PSL(2,\C)\right)/PSL(2,\C)$.}
\vspace{.1in}

The space of projective monodromy
representations with fixed traces of the $n$
generators corresponding to the punctures depends on $6g+2n-6$ complex
parameters.
These parameters are traces of certain elements of the projective
monodromy group.
The condition that monodromy is unitarizable is
that all traces are real and satisfy certain inequalities.
Thus the restriction of $p$ in Theorem 4 on the set of equations
with fixed singularities and angles is a holomorphic immersion
of a complex manifold of dimension $3g+n-3$ to a complex manifold
of dimension $6g+2n-6$. The condition of unitarizability imposes
$6g+2n-6$ {\em real} equations. The main part of Conjecture 1 is that the
set of solutions of these equations is discrete.
There are ``explicit'' expressions of the derivative of $p$ in
\cite{Hej,Earl}
but they are too complicated to obtain the necessary conclusion.
So far, the
discreteness of the set of solutions has been proved only
in the simplest cases $(g,n)=(0,4)$ and $(g,n)=(1,1)$
when we have two real equations on
one complex variable
\cite{ER}.

The problem is similar to the problem investigated by Klein and Poincar\'e
in their attempts to prove the Uniformization theorem. They tried to
show that one can choose accessory parameters so that the resulting
monodromy group is Fuchsian (this is also described by reality
of traces plus some inequalities). Eventually the Uniformization theorem
was proved by other methods. The approach of Klein and Poincare
using a Fuchsian equation has been
recently completed in the book \cite{SG}. But the proof of the
required transversality property of the monodromy map uses 
the properties of the hyperbolic metric \cite[VIII.5.3]{SG} in an essential way.

Following \cite{MP1} we denote by $M(g,\alpha_1,\ldots,\alpha_n)$
the moduli space of metrics of curvature $1$ on $S$ with
conic singularities with angles $\alpha_j$. It follows from Theorem 4,
that when none of the $\alpha_j$ are integers, and $2g+n-2>0$,
$M(g,\alpha_1,\ldots,\alpha_n)$ is a real analytic manifold of dimension
$6g-6+2n$.

A Riemannian metric defines the conformal structure on $S\backslash A$,
so we have the {\em forgetful map}
$$\Phi: M(g,\alpha_1\ldots,\alpha_n)\to M_{g,n},$$
where $M_{g,n}$ is the moduli space of $n$-punctured Riemann surfaces
of genus $g$.

To state the main results of \cite{MP1} we need
the following definitions.
$${\mathrm{Crit}}_\alpha=\left\{\|\alpha_I\|_1-\|\alpha_{I^c}\|_1+2b:
I\subset\{1,\ldots,n\}, b\in\Z_{\geq0}\right\},$$
where $\|\alpha_I\|_1=\sum_{j\in I}\alpha_j$ and $I^c=\{1,\ldots,n\}\backslash I$.
Then the {\em non-bubbling parameter} is defined by
$$NB_{g,\alpha}=d_\R\left(\chi(S\backslash A),
{\mathrm{Crit}}_\alpha\right),$$
where $d_R$ is the ordinary distance between subsets of real line.
\vspace{.1in}

\noindent
{\bf Theorem 5.} (Mondello and Panov \cite{MP1})
{\em If $NB_{g,\alpha}>0$ then the forgetful map $\Phi$ is proper.}
\vspace{.1in}

This means that under the condition $NB_{g,\alpha}>0$ the metric cannot
degenerate unless the conformal modulus degenerates. On the other hand
without this non-bubbling condition, such a degeneration is possible,
see section 10.  So to prove the main Conjecture 1 in section 2 under the
condition that $NB_{g,\alpha}>0$ it remains to prove that the set of
accessory parameters corresponding to fixed position of singularities
and fixed angles is discrete.

Another result in Mondello and Panov \cite{MP1} is that
the moduli space $M(g,\alpha_1,\ldots,\alpha_n)$ may be disconnected
for some choice of $\alpha_j$.
\vspace{.1in}

Topology of the the moduli space $M_1(\alpha_1)$ of tori with one singularity
has been described in \cite{EMP}. When $\alpha$ is not an odd integer,
denote $m=\lfloor(\alpha_1+1)/2\rfloor$. Then $M_1(\alpha_1)$ is a connected surface
of genus $\lfloor (m^2-6m+12)/12\rfloor$ with $m$ punctures.

When $\alpha_1$ is an odd integer, the metrics are coaxial, and
the set of equivalence classes has $\lceil m(m+1)/6\rceil$ connected
components each of which is an open disk.

When $\alpha_1$ is an even integer, $M_1(\alpha_1)$ has a natural
complex analytic structure (it is an algebraic curve), the forgetful
map is complex analytic (algebraic) and its degree is $\alpha_1/2$.
\vspace{.2in}

\noindent
{\bf 5. Topological degree}.
\vspace{.2in}

One can use the topological (Leray--Schauder) degree of the equation
\begin{equation}\label{eca}
\Delta\log\rho+\rho^2=2\pi\sum_{j=0}^{n-1}(\alpha_j-1)\delta_{a_j}
\end{equation}
for the density $\rho$ of the metric to obtain a lower estimate of
the number of solutions. 
\vspace{.1in}

\noindent
{\bf Theorem 6.}\cite{BdMM} {\em Let $S$ be a compact surface of genus $g>0$.
Suppose that the angles $\alpha_1,\ldots,\alpha_n$ satisfy $\alpha_j>1,\; 1\leq j\leq n$, 
$$\chi(S)+\sum_{j=1}^n(\alpha_j-1)>2\min\{\alpha_1,\ldots,\alpha_n,1\},$$
and
\begin{equation}\label{con}
\sum_{j\in I}\alpha_j-\sum_{j\not\in I}\alpha_j\neq 2k-2+n+2g,\quad k\in\Z.
\end{equation}
Then there is at least one metric on $S$ with conic singularities
at any given points with angles $\alpha_j$.}
\vspace{.1in}

For example, there is always a metric on torus with a single
singularity where the angle is $2\pi\alpha$ and $\alpha$ is not
an odd integer. On the other hand, we know from the previous section
that when $\alpha=3$ such metric may exist or not,
depending on the torus.

Condition (\ref{con}) coincides with the condition $NB_{g,\alpha}>0$
in Theorem 5. Chen and Lin \cite{CL} actually computed the degree
and included the case $g=0$.
They define the generating function
$$g(x)=(1+x+x^2+\ldots)^{-\chi(S)+n}\prod_{j=1}^n\left(1-x^{\alpha_j}\right),$$
where $\chi(S)=2-2g$. Suppose that
$$g(x)=1+b_1x^{n_1}+b_2x^{n_2}+\ldots+b_kx^{n_k}+\ldots,$$
(this is the definition of $b_k$).
\vspace{.1in}

\noindent
{\bf Theorem 7.} {\em Let $d$ be the Leray--Schauder degree of (\ref{eca}).
Define $k$ by the inequalities
$$2n_k<\chi(S)+\sum_j(\alpha_j-1)<2n_{k+1},$$
(this is well defined if (\ref{con}) holds).
Then 
$$d=\sum_{j=0}^k b_j.$$
}
\vspace{.1in}
The lower estimate of the number
of metrics which follows from degree computation is sometimes best possible.
\vspace{.1in}

We will see in Section~9 that the lower estimate of the number of metrics
which follows from degree computation is sometimes best possible.
This is not surprising since sometimes the forgetful map is complex analytic.
\vspace{.2in}

\noindent
{\bf 6. All angles on the sphere are less than 1.}
\vspace{.2in}

\noindent
{\bf Theorem 8.} (Luo and Tian, \cite{LT}) {\em A metric of curvature $1$
on the sphere
with prescribed singularities $a_j$ and angles $\alpha_j<1$ exists
if an only if 
\begin{equation}\label{lt}
0<2+\sum_{j=1}^n(\alpha_j-1)\leq 2 \min_j\{\alpha_j\}.
\end{equation}
Such a metric is unique.}
\vspace{.1in}

The LHS inequality is (\ref{gb1}) and the RHS inequality is equivalent to
(\ref{cc}) for this case. Earlier Troyanov \cite{T1} proved sufficiency
of (\ref{lt}). Luo and Tian combined PDE arguments
with geometry by considering a convex polytope in $S^3$ associated with
the metric in question.
\vspace{.2in}

\noindent
{\bf 7. All angles on the sphere are integers.}
\vspace{.2in}

In this case, the monodromy is trivial, and the
developing map $f$ is a rational function. The singularities are the
critical
points of this rational function of multiplicity $\alpha_j-1$. From
the Riemann--Hurwitz relation we obtain that
$$2+\sum_{j=1}^n(\alpha_j-1)=2d,$$
where $d$ is the degree of $f$, and another restriction is $\alpha_j\leq d$
for all $j$, which can be obtained from (\ref{cc}). So the main problem is
reduced to the question: {\em how many equivalence classes of rational
functions exist with prescribed critical points of
prescribed multiplicities?}

The answer is known for generic position of critical points:
it is the Kostka number $K(\alpha_1-1,\ldots,\alpha_n-1)$ which can
be defined as follows. Consider a rectangular diagram of size
$2\times (2d-2)$, and fill it with numbers $1,2,\ldots,n,$ so that
the number $k$ is used 
$\alpha_k-1$ times, 
and so that the entries are strictly increasing in columns, and non-decreasing
in rows. The number of obtained tables is the Kostka number $K$.

This result is due to Scherbak \cite{Sch}. When $\alpha_j=2$ for all $k$,
the Kostka number is the Catalan number:
$$\frac{(2d-2)!}{d(d-1)!},$$
and in this special case the result was obtained by L. Goldberg \cite{G}.

It follows from the results in \cite{G,Sch} that there is always
at least one metric and at most $K$ equivalence classes of them.

The simplest example of non-uniqueness is obtained by considering a rational
function of degree $3$, there are two such non-equivalent functions
sharing the position of their $4$ simple critical points. There are only two
exceptional positions of critical points,
when there is only one class (see \cite{G}).

An interesting special case is when all critical points of $f$ lie on
a circle (a circle on the Riemann sphere is a set whose points are fixed
by an anti-conformal involution; this notion does not depend on the metric).
In this case, the number of classes of metrics is exactly equal to
the Kostka number, for any location of singularities on a circle \cite{EGSV}.
Moreover, each class contains a metric which is symmetric with respect
to this circle.

So in the case of integer angles our main question has a satisfactory
solution.
\vspace{.2in}

\noindent
{\bf 8. All but two angles on the sphere are integers.}
\vspace{.2in}

We mention that a metric on the sphere cannot have one
non-integer angle 
(this is seen from the monodromy consideration).

In the case of two non-integer angles (and all the rest integers)
the monodromy is co-axial, and
the developing map has the form $z^\beta g(z)$ where $g$ is a rational function.
The result of Scherbak mentioned in the previous section has been generalized
to this case in \cite{EGT1}.
\vspace{.1in}

\noindent
{\bf Theorem 9.} {\em Let $a_1,\ldots,a_n$ and $\alpha_1,\ldots,\alpha_n$
be given, where $\alpha_1,\alpha_2$ are not integers, while $\alpha_3,\ldots,
\alpha_n$ are integers. Then there is always at least one, and at most
$E(\alpha_1,\ldots,\alpha_n)<\infty$
classes of metrics of curvature $1$ on the sphere
with these singularities and angles. For generic singularities equality
holds. In the case when $a_1,a_2,\ldots,a_n$ lie on a circle, in this order,
we have an equality, and each class contains a metric which is 
symmetric with
respect to this circle.}
\vspace{.1in}

The number $E$ can be expressed in terms of Kostka numbers,
but the expression
is somewhat complicated.
\vspace{.2in}

\noindent
{\bf 9. All but three angles on the sphere are integers.}
\vspace{.2in}

First we mention an early result which
completely solves the problem for the sphere with three singularities
\cite{E1,FKKRUY}. In this case
the location of singularities is of course irrelevant, and existence
of the metric can be obtained  from the recent results \cite{MP}, \cite{E2}
described
in section 1. To this one can add that there is always a unique class
of such metrics.

In the case when only three singularities have non-integer angles, say
$\alpha_1,\alpha_2,\alpha_3$,
the Gauss-Bonnet and the closure conditions can be written as one
inequality,
\begin{equation}\label{et}
\cos^2\pi\alpha_1+\cos^2\pi\alpha_2+\cos^2\pi\alpha_3+2(-1)^\sigma
\cos\pi\alpha_1\cos\pi\alpha_2\cos\pi\alpha_3\leq 1,
\end{equation}
where
$$\sigma=\sum_{j=4}^n(\alpha_j-1).$$

\noindent
{\bf Theorem 10.} \cite{ET} {\em If $\alpha_1,\alpha_2,\alpha_3$
are not integers, while $\alpha_4,\ldots,\alpha_n$ are integers,
then the necessary and sufficient condition of existence of non-coaxial
metric with these angles is (\ref{et}) with strict inequality.
The number of classes of such metrics is at least $1$ and at most
$\alpha_4\cdot\ldots\cdot \alpha_n$. The upper bound is attained
for generic location of singularities.}
\vspace{.1in}

\noindent
{\bf Example.} Consider the angles $(1/2,1/2,1/2,m)$, where $m$
is an integer. Theorem 10 implies that there are $m$ metrics
with these angles for generic location of singularities.
These metrics can be lifted on the torus via the $2$-to-$1$ covering
ramified over the four singular points. The resulting metric
on the torus has one singularity with
angle $2m$, and we can apply Theorem 7 to compute
the degree in this case. Using the notation introduced
before Theorem 7, we obtain
$$g(x)=(1+x+\ldots)(1-x^{2m})=1+x+\ldots+x^{2m-1}-x^{2m}+\ldots,$$
so $k=m-1$ and
$$d=\sum_0^{m-1}1=m.$$
So in this case, the degree is equal to the number of metrics,
and all metrics on the torus come from the sphere via the lifting,
so they are invariant with respect to the conformal involution of the torus.
\vspace{.1in}

All results in sections 7--9 are all of the same type:
once the restriction
on the angles is satisfied, 
a metric exists with any position of singularities, the number of classes
of metrics
is always finite, and this number is constant for generic
singularities. The general reason for this is that the equation on
the accessory parameters which gives unitarizable monodromy is algebraic.
We conjecture that there are no other cases when the equation
for accessory parameters is algebraic.

In the next section we address the only case
studied so far when the equation
on the accessory parameter is transcendental,
and for this case we will see
that existence of the metric depends on the location of singularities
in an essential way.
\vspace{.2in}

\noindent
{\bf 10. Angles $(1/2,1/2,1/2,3/2)$ on the sphere, or $3$ on a torus.}
\vspace{.2in}

Let the singularities be $a_1,a_2,a_3,a_4$. Consider the $2$-sheeted
ramified covering
$\pi:\T\to S$ by a torus with critical values $a_j$. The metric pulls back
to the torus, and we obtain a metric $\rho^*|dz|$ on $T$ with one
singularity with angle $3$. If we set 
$$\rho^*(z)=\frac{1}{\sqrt{2}}e^{u(z)/2},$$ then
$u$ will satisfy
\begin{equation}\label{pde}
\Delta u+e^{2u}=8\pi\delta,
\end{equation}
where $\delta$ is the delta function at $0$. 
Metrics obtained by such pull-back are even, so
$$u(z)=u(-z).$$
Conversely, any even solution of (\ref{pde}) corresponds to a metric on
the sphere with $4$ conic singularities with angles $(1/2,1/2,1/2,3/2)$.
The metric $\rho^*$ is coaxial, so each solution of (\ref{pde}) comes
with a one-parametric family of such solutions. This family contains exactly
one even metric which corresponds to a metric on the sphere,
see \cite{LW}, \cite[Theorem 1]{EG}.

Equation (\ref{pde}) corresponds to the Lam\'e equation on the torus
$\T$ 
\begin{equation}\label{lame}
w''-\left(2\wp(z)+\lambda\right)w=0.
\end{equation}
We denote by $F=w_1/w_2$ a ratio of two linearly independent solutions,
this is the developing map of the metric on the torus, and $F=f\circ\pi$,
where $f$ is the developing map on the sphere.
\vspace{.1in}

The question now becomes: {\em how many values of accessory parameter $\lambda$
exist for which the projective monodromy of (\ref{lame}) is unitarizable?}
\vspace{.1in}

The answer depends on the parameter $\tau=\omega_1/\omega_2$ of the torus,
where we denote the fundamental periods by $2\omega_1,2\omega_2$.
We also set $\omega_3=\omega_1+\omega_2$, $e_j=\wp(\omega_j)$,
and $\eta_j$ is defined by $\zeta(z+\omega_j)=\zeta(z)+\eta_j$,
where $\zeta$ is
the Weierstrass $\zeta$-function. 
\vspace{.1in}

\noindent
{\bf Theorem 11.} {\em For a given $\tau$, there is at most one $\lambda$
for which (\ref{lame}) has unitarizable projective monodromy. The region
in the $\tau$-plane for which such $\lambda$ exists is explicitly described
by the inequalities:
\begin{equation}\label{ineq}
\Ima\left(\frac{2\pi i}{e_j\omega_1^2+\eta_1\omega_1}-\tau\right)<0,\quad j=1,2,3.
\end{equation}
}

The region defined by (\ref{ineq}) is shaded in Fig.~1.
\vspace{.1in}

{\em Proof.} Hermite found an explicit formula for the general solution
of (\ref{lame}) see, for example \cite[Ch. II, 59]{Halphen}. Let $a$
be a solution of $\wp(a)=\lambda$. Then
$$w_{1,2}=e^{\mp z\zeta(a)}\frac{\sigma(z\pm a)}{\sigma(z)}$$
is a fundamental set of solutions of (\ref{lame}). So their ratio is
$$F(z)=e^{2z\zeta(a)}\frac{\sigma(z-a)}{\sigma(z+a)}.$$
To determine the projective monodromy we compute $F(z+2\omega)$ where
$2\omega$ is a fundamental period:
$$F(z+2\omega)=e^{4\omega\zeta(a)-4\eta a}F(z),$$
where we used the formula
$$\sigma(z+2\omega)=-e^{2\eta(z+\omega)}\sigma(z).$$
So the monodromy is unitarizable if and only if both expressions
$$\omega_1\zeta(a)-\eta_1a\quad\mbox{and}\quad\omega_2\zeta(a)-\eta_2a\quad
\mbox{are pure imaginary}.$$
This means that two equations with respect to $a$ and $\zeta=\zeta(a)$
hold:
\begin{equation}\label{re1}
\omega_1\zeta+\overline{\omega_1}\overline{\zeta}-\eta_1a-\overline{\eta_1}\overline{a}=0,
\end{equation}
and
\begin{equation}\label{re2}
\omega_2\zeta+\overline{\omega_2}\overline{\zeta}-\eta_2z-\overline{\eta_2}
\overline{a}=0.
\end{equation}
Eliminating $\overline{\zeta}$ we obtain
one linear equation of the form
\begin{equation}\label{main}
Aa+B\overline{a}+\zeta(a)=0,
\end{equation}
where
\begin{equation}\label{AB}
A=\frac{\pi}{4\omega_1^2\Ima\tau}-\frac{\eta_1}{\omega_1},\quad
B=-\frac{\pi}{2|\omega_1|^2\Ima\tau}.
\end{equation}
These constants are uniquely defined by the condition that
our equations (\ref{re1}), (\ref{re2}) are
invariant with respect to the substitution
$$(a,\zeta)\mapsto (a+2\omega_k,\zeta+2\eta_k).$$
Equation (\ref{main}) must be solved with respect
to $a$.
It was proved in \cite{BE} that besides the three
trivial solutions $a=\omega_k,\;
1\leq k\leq 3$, equation (\ref{main}) has either none or two
solutions of the form $\pm a$. Trivial solutions do not
define linearly independent $w_1$ and $w_2$ (function $F$ 
is constant).
The two non-trivial
solutions $\pm a$ of (\ref{main}), when exist, define the same
$\lambda=\wp(a)$. This proves that at most one such $\lambda$ exists
for any torus. The region $D$ in the space of
tori in which such $\lambda$ exists
is described in \cite{BE}. 
The explicit description of $D$ is the following for all $j\in\{1,2,3\},$
$$e_j\omega_1^2+\eta_1\omega_1\neq 0\quad\mbox{and}\quad\Ima\left(\frac{\pi i}{
e_j\omega_1^2+\eta_1\omega_1}-2\tau\right)<0.$$
Theorem 11 follows.

{\em Department of Mathematics,

Purdue University,

West Lafayette IN 47907 USA

www.math.purdue.edu$/\,\tilde{}$eremenko
}
\end{document}